\documentclass[11pt]{article}
\usepackage{amscd,amsfonts,amssymb,multicol,mathtext}
\usepackage[tbtags]{amsmath}
\begin{document}

\title{On the Malgrange isomonodromic deformations of \\ non-resonant
meromorphic $(2\times2)$-connections}

\author{Yuliya P.~Bibilo, Renat R.~Gontsov}

\date{}
\maketitle

\begin{abstract}
We study the tau-function and theta-divisor of an isomonodromic
family of linear differential $(2\times2)$-systems with
non-resonant irregular singularities. In some particular case the
estimates for pole orders of the coefficient matrices of the
family are applied.
\end{abstract}

\section{Introduction}

Consider a {\it meromorphic linear $(2\times2)$-system} on the
Riemann sphere $\overline{\mathbb C}$, i.~e., a system of two
linear ordinary differential equations with singularities
$a_1^0,\ldots, a_n^0\in{\mathbb C}$ and possibly $\infty$. By a
conformal mapping one can always arrange that all the
singularities are in the complex plane only. This means that one
can reduce the system to the form
\begin{eqnarray}\label{linsyst}
\frac{dy}{dz}=B(z)\,y,\qquad
B(z)=\sum_{i=1}^n\sum_{j=1}^{r_i+1}\frac{B^0_{ij}}{(z-a^0_i)^j},
\end{eqnarray}
where $y(z)\in{\mathbb C}^2$, $B_{ij}^0$ are $(2\times
2)$-matrices and $\sum_{i=1}^nB_{i1}^0=0$, to ensure that $\infty$
is not a singular point.

The non-negative integers $r_1,\ldots,r_n$ are called the {\it
Poincar\'e ranks} of the singularities $a_1^0,\ldots,a_n^0$
respectively. One can assume that the Poincar\'e ranks
$r_1,\ldots,r_m$ are positive and $r_{m+1}=\ldots=r_n=0$ (that is,
the singular points $a_{m+1}^0,\ldots,a_n^0$ are {\it Fuchsian})
for some $0\leqslant m\leqslant n$.

We consider the {\it non-resonant} case. This means that the
leading term $B^0_{i,r_i+1}$ of each non-Fuchsian singularity
$a_i^0$, $i=1,\ldots,m$, has two distinct eigenvalues. In that
case the singular points $a_1^0,\ldots,a_m^0$ are {\it irregular}.

The system (\ref{linsyst}) can be thought of as a meromorphic
connection $\nabla^0$ (more precisely, as an equation for
horizontal sections with respect to this connection) on a
holomorphically trivial vector bundle $E^0$ of rank $2$ over
$\overline{\mathbb C}$. As known (see \cite[\S21]{IY}), in a
neighbourhood of each (non-resonant) irregular singularity $a_i^0$
the local connection form $\omega^0=B(z)dz$ of $\nabla^0$ is
formally equivalent to the 1-form
$$
\omega_{\Lambda_i^0}=\sum_{j=1}^{r_i+1}\frac{\Lambda^0_{ij}}{(z-a_i^0)^j}\,dz,
$$
where $\Lambda^0_{i1},\ldots,\Lambda^0_{i,r_i+1}$ are diagonal
matrices and the leading term $\Lambda^0_{i,r_i+1}$ is conjugated
to $B^0_{i,r_i+1}$. This means that there is an invertible matrix
formal Taylor series $\widehat F(z)$ in $(z-a_i^0)$ such that the
transformation $\tilde y=\widehat F^{-1}(z)y$ transforms the
1-form $\omega^0$ into $\omega_{\Lambda_i^0}$:
$$
\omega_{\Lambda_i^0}=\widehat F^{-1}\omega^0\widehat F-\widehat
F^{-1}(d\widehat F).
$$

One should note that formally equivalent systems in a
neighbourhood $O_{a^0_i}$ of an irregular singularity $a^0_i$ are
not necessary holomorphically or meromorphically equivalent. The
system (\ref{linsyst}) has in $O_{a^0_i}$ a formal fundamental
matrix of the form
\begin{eqnarray}\label{formfund}
\widehat Y(z)=\widehat F(z)(z-a^0_i)^{\Lambda^0_{i1}}e^{Q(z)},
\qquad Q(z)=\sum_{j=1}^{r_i}\frac{\Lambda^0_{i,j+1}}{-j}\,
(z-a^0_i)^{-j}.
\end{eqnarray}
One can cover $O_{a^0_i}$ by a set of sufficiently small sectors
$S_1,\ldots,S_N$ with vertices at $a_i^0$ such that in each $S_k$
there exists a unique fundamental matrix
$Y_k(z)=F_k(z)(z-a^0_i)^{\Lambda^0_{i1}}e^{Q(z)}$ of the system
with $F_k(z)$ having $\widehat F(z)$ as an asymptotic series in
$S_k$ (see \cite[\S21]{IY}). In every intersection $S_k\cap
S_{k+1}$ the fundamental matrices $Y_k(z)$, $Y_{k+1}(z)$ are
connected by a constant matrix $C_k$: $Y_{k+1}(z)=Y_k(z)C_k$,
which is called a {\it Stokes' matrix}. If $a^0_i$ is a
non-resonant singularity, then two formally equivalent systems are
holomorphically equivalent in $O_{a^0_i}$ if and only if they have
the same sets of Stokes' matrices (see \cite[\S21]{IY} again).

Further we will focus on deformations of the system
(\ref{linsyst}) (of the pair ($E^0$, $\nabla^0$)) that allow the
local formal equivalence class
$$
\omega_{\Lambda_i}=\sum_{j=2}^{r_i+1}\frac{\Lambda_{ij}}{(z-a_i)^j}\,dz+
\frac{\Lambda^0_{i1}}{z-a_i}\,dz, \qquad i=1,\ldots,m,
$$
to vary in the sense that the diagonal matrices
$\Lambda_{i2},\ldots,\Lambda_{i,r_i+1}$ vary in a neighbourhood of
$\Lambda^0_{i2},\ldots,\Lambda^0_{i,r_i+1}$ with $\Lambda^0_{i1}$
held fixed. Thus for the set $\Lambda_i=\{\Lambda_{i2},\ldots,
\Lambda_{i,r_i+1}\}$ of $r_i$ diagonal matrices we denote by
$\nabla_{\Lambda_i}$ the meromorphic connection on the
holomorphically trivial vector bundle of rank $2$ over $O_{a_i}$
whose 1-form is $\omega_{\Lambda_i}$. To describe the required
deformations in more details let us begin with a deformation
space.

For $k\in\mathbb N$ we denote by $Z^k$ the subset of the space
${\mathbb C}^k$ whose points have pairwise distinct coordinates.
Then $Z^n$ will be the space of pole locations and
$$
\mathcal{C}_i=\underbrace{{\mathbb C}^2\times\ldots\times{\mathbb
C}^2}_{r_i-1}\times Z^2, \qquad i=1,\ldots,m,
$$
will be the space of local formal equivalence classes at the pole
$a_i$ (any class is determined by $r_i-1$ diagonal matrices
$\Lambda_{i2},\ldots,\Lambda_{i,r_i}$ and a diagonal matrix
$\Lambda_{i,r_i+1}$ whose eigenvalues are pairwise distinct).
Define the deformation space $\mathcal{D}$ as the universal cover
$$
\mathcal{D}=\widetilde{Z^n}\times\widetilde{\mathcal{ C}}_1\times\ldots
\times\widetilde{ \mathcal{C}}_m
$$
of the Cartesian product ${Z^n}\times \mathcal{C}_1\times\ldots
\times \mathcal{C}_m$.

One has the standard projections
\begin{eqnarray*}
a=(a_1,\ldots,a_n): \mathcal{D} & \rightarrow & Z^n, \\
\Lambda_i=(\Lambda_{i2},\ldots,\Lambda_{i,r_i+1}): \mathcal{D} &
\rightarrow & \mathcal{C}_i, \qquad i=1,\ldots,m.
\end{eqnarray*}
For $t\in\mathcal{D}$ we denote by $a_i(t)$ the $i$-th coordinate of
the image of $t$ under the first projection and by $\Lambda_i(t)$
the image of $t$ under the second one. Denote then by $t^0$ the
base point of the deformation space $\mathcal{D}$ corresponding to the
system (\ref{linsyst}) (to the initial connection $\nabla^0$),
i.~e., $a(t^0)=(a_1^0,\ldots,a_n^0)$, $\Lambda_i(t^0)=
(\Lambda^0_{i2},\ldots,\Lambda^0_{i,r_i+1})$. Consider also the
singular hypersurfaces
$$
X_i=\{(z,t)\in\overline{\mathbb C}\times \mathcal{D}\mid
z=a_i(t)\}\subset\overline{\mathbb C}\times \mathcal{D}, \qquad
i=1,\ldots,n.
$$

Now for $i=1,\ldots,m$, consider the fibre bundle
$\mathcal{M}_i\rightarrow \mathcal{C}_i$, whose fiber over each
point $\Lambda_i\in \mathcal{C}_i$ is the moduli space of local
{\it holomorphic} equivalence classes of connections that are all
formally equivalent to the connection $\nabla_{\Lambda_i}$. A
point of this fiber (a holomorphic equivalence class of
connections) is determined by a corresponding set of Stokes'
matrices. Let $\sigma_i^0\in \mathcal{M}_i$ denote the holomorphic
equivalence class of the connection
$\nabla^0|_{O_{a_i^0}}\sim\nabla_{\Lambda_i^0}$ and let $\sigma_i$
denote the unique horizontal section of the fibre bundle
$\mathcal{M}_i\rightarrow \mathcal{C}_i$ such that
$\sigma_i(\Lambda_i^0)=\sigma_i^0$.

Due to B.\,Malgrange \cite[Th. 3.1]{Ma} (see also \cite[Th.
2.9]{Pa}) the following statement holds.
\medskip

{\bf Theorem 1.} {\it There exists a} {\rm unique}\footnote{Under
some additional assumption we discuss later on.} {\it
isomonodromic deformation $(E,\nabla)$ of the pair
$(E^0,\nabla^0)$, that is, the rank $2$ holomorphic vector bundle
$E$ over $\overline{\mathbb C} \times\mathcal{D}$ and integrable
meromorphic connection $\nabla$ on $E$ with a} simple {\it type
$r_i$ singularity\footnote{That is, near $X_i$ the local
connection 1-form $\Omega$ of $\nabla$ looks like
$$
\Omega=\frac{B_i(z,t)}{(z-a_i(t))^{r_i+1}}\,d(z-a_i(t))+
\sum_k\frac{C_{ik}(z,t)}{(z-a_i(t))^{r_i}}\,dt_k,
$$
where the matrices $B_i$, $C_{ik}$ are holomorphic and (for
$i=1,\ldots,m$) the eigenvalues of $B_i(a_i(t),t)$ are pairwise
distinct.} along $X_i$, $i=1,\ldots,n$, satisfying the following
properties:
\begin{itemize}
\item the restriction of $(E,\nabla)$ to $\overline{\mathbb
C}\times\{t^0\}$ is isomorphic to $(E^0,\nabla^0)$;

\item for any $t\in\cal D$ the restriction of $\nabla$ to
$\overline{\mathbb C}\times\{t\}$ is formally equivalent to the
local connection $\nabla_{\Lambda_i(t)}$ near $z=a_i(t)$,
$i=1,\ldots,m$, and belongs to the local holomorphic equivalence
class $\sigma_i(\Lambda_i(t))\in\mathcal{M}_i$.
\end{itemize}
}

The deformation described above will be referred to as the {\it
Malgrange isomonodromic deformation} of the pair $(E^0,\nabla^0)$.

According to the Malgrange--Helminck--Palmer theorem (see
\cite[\S3]{Pa} or \cite[\S3]{Ma}) the set
$$
\Theta=\{t\in \mathcal{D}\mid E|_{\overline{\mathbb C}\times\{t\}}
\hbox{ is non-trivial }\}
$$
is either empty or $\Theta\subset\mathcal{D}$ is an analytic
subset of codimension one (which is usually called the {\it
Malgrange $\Theta$-divisor}). If the latter holds, there exists a
function $\tau$ (called the $\tau$-{\it function} of the
isomonodromic deformation) holomorphic on the whole space $\cal D$
whose zero set coincides with $\Theta$.

Thus the Malgrange isomonodromic deformation of the pair
$(E^0,\nabla^0)$ determines an isomonodromic deformation
\begin{eqnarray}\label{isomfam}
\frac{dy}{dz}=\Bigl(\sum_{i=1}^n\sum_{j=1}^{r_i+1}
\frac{B_{ij}(t)}{(z-a_i(t))^j}\Bigr) y, \qquad B_{ij}(t^0)
=B^0_{ij},
\end{eqnarray}
of the system (\ref{linsyst}) for $t\in D(t^0)$, where $D(t^0)$ is
a neighbourhood of the point $t^0$ in the space $\mathcal{D}$. The
matrix functions $B_{ij}(t)$, holomorphic in $D(t^0)$, can be
extended meromorphically to the whole space $\mathcal{D}$ and have
$\Theta$ as a polar locus.

Recall that for a Fuchsian system (the case of $m=0$)
\begin{eqnarray}\label{linsystfuchs}
\frac{dy}{dz}=\Bigl(\sum_{i=1}^n\frac{B^0_i}{z-a^0_i}\Bigr)y
\end{eqnarray}
the best known isomonodromic deformation has been described by
L.\,Schlesinger \cite{Sch}, \cite{Sch2}. Starting from the initial
conditions $B_i(a^0)=B_i^0$, $a^0=(a^0_1,\ldots,a^0_n)$, the
residue matrices $B_i(a)$ vary satisfying the {\it Schlesinger
equation}
$$
dB_i(a)=-\sum_{j=1,j\ne i}^n\frac{[B_i(a),B_j(a)]}{a_i-a_j}\,
d(a_i-a_j), \qquad i=1,\ldots,n,
$$
and they are extended as meromorphic matrix functions to the
deformation space $\widetilde{Z^n}$ from a neighbourhood $D(a^0)$
of the initial point $a^0$.

A.\,A.\,Bolibruch \cite{Bo2} has obtained the following result
concerning pole orders of the matrices $B_i(a)$.
\medskip

{\bf Theorem 2.} {\it Let the monodromy of the $(2\times2)$-system
$(\ref{linsystfuchs})$ be} irreducible {\it and let $a^*\in\Theta$
be a point of the $\Theta$-divisor such that the restriction
$E|_{\overline{\mathbb C}\times\{a^*\}}$ is of the form
$$
E|_{\overline{\mathbb C}\times\{a^*\}}\cong{\cal O}(-1)\oplus{\cal
O}(1).
$$
Then in a neighbourhood $D(a^*)$ of $a^*$ the $\Theta$-divisor is
an analytic} submanifold {\it and the matrix functions $B_i(a)$
have poles of at most} second {\it order along
$D(a^*)\cap\Theta$.}

The latter means that $\tau^2(a)B_i(a)$ are holomorphic matrix
functions in $D(a^*)$. The proof of Theorem 2 (formulated in a
more general setting) also contains in \cite{GV}.
\medskip

Adapting Bolibrukh's ideas to the case of linear systems with
irregular singularities we propose a local description of the
$\Theta$-divisor of the Malgrange isomonodromic deformation and
generalization of Theorem 2 when the initial system has at most
two irregular singularities and their Poincar\'e ranks are equal
to $1$ (Theorem 3).

\section{Holomorphic vector bundles and the Riemann--Hilbert problem for irregular systems}

The fact $t^*\in\Theta$ means that the restriction
$E|_{\overline{\mathbb C}\times\{t^*\}}$ of the holomorphic vector
bundle $E$ described in Theorem 1 is not holomorphically trivial.
This restriction belongs to the family $\cal F$ of holomorphic
vector bundles over the Riemann sphere endowed with meromorphic
connections which occurs in the investigation of the corresponding
Riemann--Hilbert problem. The latter is the question on existence
of a global meromorphic linear system with the singular points
$a_1^*=a_1(t^*),\ldots,a_n^*=a_n(t^*)$ of Poincar\'e ranks
$r_1,\ldots,r_n$ respectively that

1) has the same monodromy as the initial one and

2) is meromorphically equivalent to the local system
\begin{eqnarray}\label{locsyst}
dy=\omega_i^*y
\end{eqnarray}
determined by the local holomorphic equivalence class
$\sigma_i(\Lambda_i(t^*))$ near each irregular singular point
$a_i^*$.

The Riemann--Hilbert problem under consideration has a positive
answer (it is sufficient one of the irregular singularities to be
non-resonant for positive solution in the two-dimensional case,
see \cite{Bol_Malek_Mit}). This means there is a holomorphic
vector bundle (not $E|_{\overline{\mathbb C}\times\{t^*\}}$) in
the family $\cal F$ that is holomorphically trivial. Thus we are
coming to the point where it is naturally to recall briefly the
construction of the family $\cal F$ (see details in
\cite{Bol_Malek_Mit}).

By the monodromy representation (generated by the monodromy
matrices $G_1,\ldots,G_n$) of the initial system (\ref{linsyst})
one constructs over the punctured Riemann sphere
$\overline{\mathbb C}\setminus\{a_1^*,\ldots,a_n^*\}$ a
holomorphic vector bundle $\widetilde F$ of rank $2$ with a
holomorphic connection $\widetilde\nabla$ having the prescribed
monodromy. This bundle is defined by a set $\{U_{\alpha}\}$ of
sufficiently small discs covering $\overline{\mathbb C}
\setminus\{a_1^*,\ldots,a_n^*\}$ and set $\{g_{\alpha\beta}\}$ of
constant matrices defining a gluing cocycle. A connection
$\widetilde\nabla$ is defined by a set $\{\omega_{\alpha}\}$ of
matrix differential 1-forms $\omega_{\alpha}\equiv0$. So in the
intersections $U_{\alpha}\cap U_{\beta}\ne\varnothing$ the gluing
conditions
$$
\omega_{\alpha}=(dg_{\alpha\beta})g_{\alpha\beta}^{-1}+
g_{\alpha\beta}\omega_{\beta}g_{\alpha\beta}^{-1}
$$
hold.

Further one extends the pair $(\widetilde F,\widetilde\nabla)$ to
the whole Riemann sphere. In neighbourhoods $O_{a_i^*}$ of the
irregular singular points $a_i^*$, $i=1,\ldots,m$, the extension
of $\widetilde\nabla$ is determined by the corresponding local
matrix differential 1-forms $\omega_i^*$ of the coefficients of
the systems (\ref{locsyst}), while in neighbourhoods $O_{a_i^*}$
of the Fuchsian singular points $a_i^*$, $i=m+1,\ldots,n$, the
extension of $\widetilde\nabla$ is determined by the matrix
differential 1-forms $E_idz/(z-a_i^*)$. Here $E_i=
1/(2\pi\sqrt{-1})\ln G_i$ is a normalized logarithm of the
monodromy matrix $G_i$ and its branch is chosen so that the
eigenvalues $\rho_i^k$ of $E_i$ satisfy the condition
\begin{eqnarray}\label{cond}
0\leqslant{\rm Re}\,\rho_i^k<1.
\end{eqnarray}
This is the so-called {\it canonical} extension $(\widetilde
F^0,\widetilde\nabla^0)$ of the pair $(\widetilde
F,\widetilde\nabla)$ in the sense of Malgrange \cite{Ma2} (and
Deligne \cite{De}, for the Fuchsian case).

Finally, consider a formal fundamental matrix (see
(\ref{formfund}))
$$
\widehat Y_i(z)=\widehat F_i(z)(z-a_i^*)^{\Lambda^0_{i1}}
e^{Q_i(z)}, $$
$$Q_i(z)=\sum_{j=1}^{r_i}
\frac{\Lambda^*_{i,j+1}}{-j}\,(z-a_i^*)^{-j}, \quad
\Lambda^*_{i,j+1} =\Lambda_{i,j+1}(t^*),
$$
of each local irregular system (\ref{locsyst}), $i=1,\ldots,m$,
and write it in the form
\begin{eqnarray}\label{formfund2}
\widehat Y_i(z)=\widehat
F_i(z)(z-a_i^*)^{A^0_i}(z-a_i^*)^{\widehat E_i}e^{Q_i(z)}, \qquad
A^0_i=[{\rm Re}\,\Lambda^0_{i1}].
\end{eqnarray}
The diagonal elements of the integer-valued matrix $A^0_i$ are
referred to as the {\it formal valuations} of the system. As
follows, the diagonal elements $\rho_i^k$ of the matrix $\widehat
E_i$ satisfy the condition (\ref{cond}). By an analogue of
Sauvage's lemma (see \cite[L. 11.2]{Ha}) for formal matrix series,
for any diagonal integer-valued matrix $A_i$ there exists a matrix
$\Gamma'_i(z)$ meromorphically invertible in $O_{a_i^*}$ such that
\begin{eqnarray}\label{sauvage}
\Gamma'_i(z)\widehat F_i(z)(z-a_i^*)^{A^0_i-A_i}=
(z-a_i^*)^{\widetilde A_i}\widehat H_i(z),
\end{eqnarray}
where $\widetilde A_i$ is a diagonal integer-valued matrix and
$\widehat H_i(z)$ is an invertible formal (matrix) Taylor series
in $z-a_i^*$.

Now one constructs the family $\cal F$ of extensions of the pair
$(\widetilde F,\widetilde\nabla)$ replacing the form $\omega_i^*$
in the construction of $(\widetilde F^0,\widetilde\nabla^0)$ by
the form
$$
\omega^{A_i}=(d\Gamma_i)\Gamma_i^{-1}+\Gamma_i\omega^*_i\Gamma_i^{-1},
\quad \Gamma_i(z)=(z-a_i^*)^{-\widetilde A_i}\Gamma'_i(z), \quad
i=1,\ldots,m,
$$
and the form $E_idz/(z-a^*_i)$ by the form
$$
\omega^{A_i}=(d\Gamma_i)\Gamma_i^{-1}+\Gamma_i\frac{E_idz}{z-a^*_i}\Gamma_i^{-1},
\quad \Gamma_i(z)=(z-a_i^*)^{A_i}S_i, \quad i=m+1,\ldots,n,
$$
where $A_i={\rm diag}(d_i^1,d_i^2)$ is a diagonal integer-valued
matrix whose diagonal elements satisfy the condition
$d_i^1\geqslant d_i^2$, and $S_i$ is a non-singular matrix
reducing the matrix $E_i$ to an upper-triangular form
$E'_i=S_iE_iS_i^{-1}$. As follows from (\ref{formfund2}),
(\ref{sauvage}), a formal fundamental matrix of the local
irregular system $dy=\omega^{A_i}y$, $i=1,\ldots,m$, is of the
form
\begin{eqnarray}\label{formfund3}
\widehat Y'_i(z)=\Gamma_i(z)\widehat{Y}_i(z)=\widehat H_i(z)(z-a_i^*)^{A_i}(z-a_i^*)^{\widehat
E_i}e^{Q_i(z)}.
\end{eqnarray}
Its singular point $z=a_i^*$ is of Poincar\'e rank $r_i$ again. At
the same time, the local system $dy=\omega^{A_i}y$,
$i=m+1,\ldots,n$, is Fuchsian:
$$
\omega^{A_i}=\Bigl(\frac{A_i}{z-a^*_i}+(z-a^*_i)^{A_i}
\frac{E'_i}{z-a^*_i}(z-a^*_i)^{-A_i}\Bigr)dz.
$$

Let us call the matrices $A_1,\ldots,A_n$, $S_{m+1},\ldots,S_n$
involved in the construction above, the {\it admissible} matrices.
Thus the family $\cal F$ consists of the pairs $(F^{A,S},
\nabla^{A,S})$ obtained by all sets $(A,S)=\{A_1,\ldots,A_n,
S_{m+1},\ldots,S_n\}$ of admissible matrices. Though the matrices
$\Gamma'_1(z),\ldots,\Gamma'_m(z)$ (see (\ref{sauvage})) are also
involved in the construction of the pair $(F^{A,S},
\nabla^{A,S})$, one should note that in our (non-resonant) case
the bundle $F^{A,S}$ does not depend on them (for a fixed
$(A,S)$).

Now the restriction $(E,\nabla)|_{\overline{\mathbb C}\times
\{t^*\}}$ can be thought of as an element of the family $\cal F$:
$$
(E,\nabla)|_{\overline{\mathbb C}\times\{t^*\}}
\cong(F^{A^0,S^0},\nabla^{A^0,S^0}), $$
$$
A^0=\{A^0_1,\ldots,A^0_n\}, \quad S^0=\{S^0_{m+1},\ldots,S^0_n\},
$$
where the matrices $A^0_1,\ldots,A^0_m$ are defined in
(\ref{formfund2}), and the sets of the (admissible) matrices
$A^0_{m+1},\ldots,A^0_n$ and $S^0_{m+1},\ldots,S^0_n$ come from
the {\it Levelt} decompositions \cite{Le} of a fundamental matrix
$Y(z)$ of the initial system (\ref{linsyst}) at the corresponding
Fuchsian singularities $a^0_{m+1},\ldots,a^0_n$:
$$
Y(z)=U_i(z)(z-a^0_i)^{A^0_i}S^0_i(z-a^0_i)^{E_i},\qquad
i=m+1,\ldots,n,
$$
where the matrix $U_i(z)$ is holomorphically invertible at the
point $a^0_i$. The matrices $A^0_{m+1},\ldots,A^0_n$ are preserved
along the deformation (see \cite{Bo3}). And one requires the
matrices $S^0_{m+1},\ldots,S^0_n$ to be also preserved, to ensure
that the Malgrange deformation is a unique isomonodromic
deformation of the pair $(E^0,\nabla^0)$ (see Theorem 1).

\section{Theorem on $\Theta$-divisor}

Now let us consider a linear meromorphic $(2\times2)$-system with
$n$ singular points such that $m\leqslant2$ of them are irregular
and their Poincar\'e ranks are equal to $1$, i.~e., the system of
the form (\ref{linsyst}), where $r_{1,2}\leqslant1$,
$r_3=\ldots=r_n=0$:
\begin{eqnarray}\label{22syst}
\frac{dy}{dz}=\biggl(\frac{B^0_{12}}{(z-a^0_1)^2}+\frac{B^0_{22}}{(z-a^0_2)^2}
+\sum_{i=1}^n\frac{B^0_{i1}}{z-a^0_i}\biggr)y.
\end{eqnarray}
The $\Theta$-divisor and the coefficient matrices $B_{ij}(t)$ of
the Malgrange isomonodromic deformation (\ref{isomfam}) of such
system possess the following properties.
\medskip

{\bf Theorem 3.} {\it Let the monodromy representation of the
$(2\times2)$-system $(\ref{22syst})$ be} irreducible {\it and let
$t^*\in\Theta$ be a point of the $\Theta$-divisor such that
$$
E|_{\overline{\mathbb C}\times\{t^*\}}\cong
\mathcal{O}(-1)\oplus\mathcal{O}(1).
$$
Then in a neighbourhood $D(t^*)$ of $t^*$ the $\Theta$-divisor is
an analytic} submanifold {\it and the matrix functions $B_{ij}(t)$
have poles of at most} second {\it order along $D(t^*)\cap
\Theta$.}
\medskip

Before proving this theorem let us recall a calculation algorithm
for the local $\tau$-function of the Malgrange isomonodromic
deformation $(E,\nabla)$ of the system (\ref{22syst}).

Consider a point $t^*\in\Theta$. Though the corresponding pair
$(E,\nabla)|_{\overline{\mathbb C}\times\{t^*\}}\cong
(F^{A^0,S^0},\nabla^{A^0,S^0})$ is such that the bundle
$$
F^{A^0,S^0}\cong{\cal O}(-1)\oplus{\cal O}(1)
$$
is not holomorphically trivial, one can construct an auxiliary
linear meromorphic system
\begin{eqnarray}\label{linsyst_aux}
\frac{dy}{dz}=\biggl(\frac{B^{*0}_{12}}{(z-a^*_1)^2}+\frac{B^{*0}_{22}}{(z-a^*_2)^2}
+\sum_{i=1}^n\frac{B^{*0}_{i1}}{z-a_i^*}\biggr)y,
\end{eqnarray}
with irregular non-resonant singular points $a^*_1=a_1(t^*)$,
$a^*_2=a_2(t^*)$ of Poincar\'e rank $1$ and Fuchsian singular
points $a^*_3=a_3(t^*),\ldots,a^*_n=a_n(t^*)$. This system is
holomorphically equivalent to the local systems determined by the
connection $\nabla^{A^0,S^0}$ in neighbourhoods of the
corresponding singular points, but it has an {\it apparent}
Fuchsian singularity at the infinity (i.~e., the monodromy at this
point is trivial). Its fundamental matrix is of the form
$Y^*(z)=U(z)z^K$ near the infinity, where
\begin{eqnarray}\label{seriesU}
U(z)=I+U_1\frac1z+U_2\frac1{z^2}+\ldots, \qquad K={\rm diag}
(-1,1).
\end{eqnarray}
Therefore, the residue matrix at the infinity is equal to $-K$,
and $\sum_{i=1}^n B^{*0}_{i1}=K$ (existence of such a system in
the Fuchsian case is explained, for example, in the proof of
Theorem 2 from \cite{VG}; an explanation here is the same).

The columns of the fundamental matrix $Y^*(z)$ of the system
$(\ref{linsyst_aux})$ under $\mathbb{C}$ determine a basis of
sections of the bundle $F^{A^0,S^0}$ horizontal with respect to
$\nabla^{A^0,S^0}$. Consider a matrix $V(z)$ holomorphically
invertible in a neighbourhood $O_{\infty}$ of the infinity whose
columns determine this basis under $O_{\infty}$. Then the quotient
$Y^*(z)V^{-1}(z)=g_{0 \infty}$ is a cocycle of the bundle
$F^{A^0,S^0}$, which is $z^K$, i.e.
\begin{eqnarray}\label{UV}
U(z)z^K=z^KV(z).
\end{eqnarray}

Let us include the auxiliary system $(\ref{linsyst_aux})$ into the
Malgrange isomonodromic family
\begin{eqnarray}\label{isom_fam_aux}
\frac{dy}{dz}=\biggl(\frac{B^*_{12}(t)}{(z-a_1(t))^2}+
\frac{B^*_{22}(t)}{(z-a_2(t))^2}+\sum_{i=1}^n\frac{B^*_{i1}(t)}
{z-a_i(t)}\biggr)y, \quad B^*_{ij}(t^*)=B^{*0}_{ij}.
\end{eqnarray}
An appropriate matrix meromorphic differential 1-form determining
this family (see \cite[Ch.4, \S1]{FIKN}) has the form
\begin{eqnarray}\label{def_form_aux}
\omega=\sum_{i=1}^2\frac{B^*_{i2}(t)}{(z-a_i(t))^2}\,d(z-a_i(t))+
\sum_{i=1}^n\frac{B^*_{i1}(t)}{z-a_i(t)}\,
d(z-a_i(t))+\;(d\Lambda)\mbox{-part}.
\end{eqnarray}
Observe that the equality $\sum_{i=1}^n B^*_{i1}(t)=K$ holds.
Indeed, the differential $1$-form $\omega$ satisfies the Frobenius
integrability condition, i.e., $d\omega=\omega\wedge\omega$. One
can directly check that the residue (in the sense of Leray) of
$\omega\wedge\omega$ along $\{z=\infty\}$ is equal to zero and the
residue of $d\omega$ along $\{z=\infty\}$ is equal to
$d\sum_{i=1}^n B^*_{i1}(t)$.

Let $Y(z,t)$ be the fundamental matrix of the Pfaffian system
$dy=\omega y$ of the form
\begin{eqnarray}\label{fundinf}
Y(z,t)=U(z,t) z^K, \quad U(z,t)=I+U_1(t)\frac1z+
U_2(t)\frac1{z^2}+\ldots,
\end{eqnarray}
at the infinity, and $Y(z,t^*)=Y^*(z)$ (by analogy with the
Fuchsian case \cite{Bo}).

As follows from $(\ref{def_form_aux})$,
\begin{eqnarray}\label{partial_Y}
\frac{\partial Y}{\partial a_i} Y^{-1}=-\sum_{j=1}^{r_i+1}
\frac{B^*_{ij}(t)}{(z-a_i)^j}=-\sum_{j=1}^{r_i+1}
\frac{B^*_{ij}(t)}{z^j(1-\frac{a_i}z)^j}.
\end{eqnarray}
Expanding into series the left and the right sides of
$(\ref{partial_Y})$ near the infinity, one gets
$$
\frac{\partial U_1(t)}{\partial a_i}\frac1z+o(z^{-1})=
\Bigl(-B^*_{i1}(t)\frac1z+o(z^{-1})\Bigr)
\Bigl(I+U_1(t)\frac1z+o(z^{-1})\Bigr),
$$
therefore
\begin{eqnarray}\label{u1_partial}
\frac{\partial U_1(t)}{\partial a_i}=-B^*_{i1}(t), \qquad
i=1,\ldots,n.
\end{eqnarray}

From the relation
$$
\frac{\partial Y}{\partial z}Y^{-1}=\sum_{i=1}^n\sum_{j=1}^{r_i+1}
\frac{B^*_{ij}(t)}{z^j(1-\frac{a_i}z)^j}
$$
one gets
\begin{eqnarray*}
-U_1(t)\frac1{z^2}+o(z^{-2})+\Bigl(I+U_1(t)\frac1z+o(z^{-1})\Bigr)\frac Kz=\\
=\Bigl(\frac Kz+\Bigl(\sum_{i=1}^nB^*_{i1}(t)a_i+B^*_{12}(t)+
B^*_{22}(t)\Bigr)\frac1{z^2}+o(z^{-2})\Bigr)
\Bigl(I+U_1(t)\frac1z+o(z^{-1})\Bigr).
\end{eqnarray*}
Hence
$$
-U_1+[U_1,K]=\sum_{i=1}^nB^*_{i1}(t)a_i+B^*_{12}(t)+B^*_{22}(t).
$$
Thus the upper-right element $u_1(t)$ of the matrix $U_1(t)$
coincides with the same element of the matrix
$\sum_{i=1}^nB^*_{i1}(t)a_i+B^*_{12}(t)+B^*_{22}(t)$.
\medskip

{\bf Lemma 1.} {\it The function $u_1(t)$ is not equal to zero
identically and vanishes at the point $t=t^*$.}
\medskip

{\bf Proof.} Since the matrix $U_1(t^*)$ is that from the
decomposition (\ref{seriesU}), the vanishing of $u_1(t)$ at the
point $t^*$ follows from the relation (\ref{UV}).

Now let us explain that the function $u_1(t)$ is not equal to zero
identically. We denote by $b_{ij}(t)$ the upper-right elements of
the matrices $B^*_{ij}(t)$. Then
$$
u_1(t)=b_{12}(t)+b_{22}(t)+\sum_{i=1}^n b_{i1}(t)a_i
$$
and as follows from (\ref{u1_partial}),
$$
\frac{\partial u_1(t)}{\partial a_i}=-b_{i1}(t), \qquad
i=1,\ldots,n.
$$
Arguing by contradiction, suppose that $u_1(t)\equiv 0$. Then the
following equalities should be true:
\begin{equation*}
\begin{split}
&b_{i1}(t)\equiv0, \qquad i=1,\ldots,n, \\
&b_{12}(t)+b_{22}(t)\equiv0.
\end{split}
\end{equation*}
We will show that $b_{12}(t)=b_{22}(t)\equiv0$ as well, which
contradicts irreducibility of the monodromy of the family
(\ref{isom_fam_aux}).

To use the fact that $z=\infty$ is an apparent singularity of the
family (\ref{isom_fam_aux}), let us turn to a new independent
variable $\xi=z^{-1}$ and examine the matrix differential 1-form
$B^*(z,t)dz$ of the coefficients of this family near the point
$\xi=0$:
\begin{equation*}
\begin{split}
&B^*(z,t)dz=-\frac{B^*(\xi^{-1},t)}{\xi^2}\,d\xi, \quad
-\frac{B^*(\xi^{-1},t)}{\xi^2}=-\sum_{i=1}^2\frac{B^*_{i2}(t)}{(1-a_i\xi)^2}-
\sum_{i=1}^n\frac{B^*_{i1}(t)}{\xi(1-a_i\xi)}=\\
&=\frac{-1}{\xi}\biggl(K+\sum_{i=1}^nB^*_{i1}(t)a_i\xi+\sum_{i=1}^n
B^*_{i1}(t)a^2_i\xi^2+o(\xi^2)\biggr)-\biggl(\sum_{i=1}^2B^*_{i2}(t)+
2\sum_{i=1}^2B^*_{i2}(t)a_i\xi+o(\xi)\biggr)=\\
&=\frac{-1}{\xi}K-\biggl(\sum_{i=1}^nB^*_{i1}(t)a_i+\sum_{i=1}^2B^*_{i2}(t)\biggr)-
\biggl(\sum_{i=1}^nB^*_{i1}(t)a^2_i+2\sum_{i=1}^2B^*_{i2}(t)a_i\biggr)\xi+o(\xi)=\\
&=\left(\begin{array}{r r}
1 & 0 \\
0 & -1\end{array}\right)\frac1{\xi}+
\left(\begin{array}{r r}
* & 0 \\
* & *
\end{array}\right)+
\left(\begin{array}{c c}
* & -2\sum_{i=1}^2b_{i2}(t)a_i \\
* & *
\end{array}\right)\xi+o(\xi).
\end{split}
\end{equation*}

The gauge transformation $\tilde{y}=\xi^K y$ changes the latter
matrix into a new one having the form
$$
\frac1{\xi}\left(\begin{array}{c c}
0 & -2\sum_{i=1}^2 b_{i2}(t)a_i \\
0 & 0
\end{array}\right)+O(1).
$$
The monodromy matrix of the Fuchsian singular point $\xi=0$ of the
transformed system is identity. On the other hand, both
eigenvalues of its residue matrix are zeros. Thus the monodromy
matrix is equal to the exponent of the residue matrix, i. e.,
\begin{equation*}
\exp 2\pi\sqrt{-1} \left(\begin{array}{c c}
0 & -2 \sum_{i=1}^2b_{i2}(t)a_i \\
0 & 0
\end{array}\right) = I.
\end{equation*}
Then the equality $b_{12}(t)a_1+b_{22}(t)a_2\equiv0$ holds, which
(together with the equality $b_{12}(t)+b_{22}(t)\equiv0$) implies
$b_{12}(t)=b_{22}(t)\equiv0$. {\hfill $\Box$}
\medskip

{\bf Lemma 2.} {\it The function $u_1(t)$ is a local
$\tau$-function of the Malgrange isomonodromic deformation of the
system $(\ref{22syst})$, i.~e., it locally determines the
$\Theta$-divisor near the point $t^*\in\Theta$.}
\medskip

{\bf Proof.} If $u_1(t)\ne0$, then we can consider a
holomorphically invertible (with respect to $z$) in $\mathbb C$
matrix
$$
\Gamma'_1(z,t)=\left(\begin{array}{cc} 1 & 0 \\
                       -\frac{z}{u_1(t)} & 1
                     \end{array}\right).
$$
By the construction the matrix $U'(z,t)=\Gamma'_1(z,t)U(z,t)$ is
of the form
$$
U'(z,t)=\Bigl(U'_0(t)+U'_1(t)\frac1z+\ldots\Bigr)z^{-K}, \qquad
U'_0(t)=\left(\begin{array}{cc} 0 & u_1(t) \\
                  -\frac1{u_1(t)} & \frac{f(t)}{u_1(t)}
              \end{array}\right),
$$
where $f(t)$ is a holomorphic function at the point $t=t^*$.

The gauge transformation
\begin{equation}\label{Gamma1}
y_1=\Gamma_1(z,t)y, \qquad
\Gamma_1(z,t)=U'_0(t)^{-1}\Gamma'_1(z,t),
\end{equation}
transforms the system (\ref{isom_fam_aux}) into a new one, with a
fundamental matrix
\begin{eqnarray}\label{triv}
Y^1(z,t)=\Gamma_1(z,t)Y(z,t)
\end{eqnarray}
that is holomorphically invertible at the infinity. As the columns
of the matrix $Y(z,t)$ form a basis of horizontal (with respect to
the restriction of the connection $\nabla$ on $\overline{\mathbb
C} \times\{t\}$) sections of the bundle $E_{\overline{\mathbb C}
\times\{t\}}$ over $\mathbb C$, the relation (\ref{triv}) implies
a holomorphic triviality of this bundle.

If $u_1(t)=0$, then the matrix
$$
V_{\infty}(z)=z^{-K}U(z,t)z^K=z^{-K}\left(I+\left(
\begin{array}{r r} * & 0 \\
                   * & *
\end{array}\right)\frac1z+\ldots\right)z^K
$$
is holomorphically invertible at the infinity, hence
$Y(z,t)=z^KV_{\infty}(z)$ and $E|_{\overline{\mathbb C}\times
\{t\}}\cong{\cal O}(-1)\oplus{\cal O}(1)$. {\hfill $\Box$}
\medskip

{\bf Proof of Theorem 3.} First we explain that
$du_1(t^*)\not\equiv0$. Indeed, in the opposite case the following
equalities should be true:
\begin{equation*}
\begin{split}
&b_{i1}(t^*)=0, \qquad i=1,\ldots,n, \\
&b_{12}(t^*)+b_{22}(t^*)=0.
\end{split}
\end{equation*}
Then similarly to the proof of Lemma 1 one gets the relations
$b_{12}(t^*)=b_{22}(t^*)=0$, which contradict the monodromy
irreducibility. Thus the $\Theta$-divisor of the Malgrange
isomonodromic deformation of the system (\ref{22syst}) is an
analytic submanifold in a neighbourhood $D(t^*)$ of the point
$t^*$.

Now let us estimate the pole orders of the matrices $B_{i1}(t)$,
$B_{12}(t)$, $B_{22}(t)$ along $\Theta\cap D(t^*)$. Return to the
proof of Lemma 2. The family obtained from (\ref{isom_fam_aux})
via the gauge transformation (\ref{Gamma1}), coincides with the
Malgrange isomonodromic deformation (for $t\in D(t^*)\setminus
\Theta$) of the initial system (\ref{22syst}). (Indeed, this
transformation does not change connection matrices at the Fuchsian
singular points and it also does not change holomorphic
equivalence classes of the family at the irregular singularities.)
Therefore the coefficient matrix of the Malgrange isomonodromic
deformation of the initial system (\ref{22syst}) has the form
$$
\frac{\partial\Gamma_1}{\partial z}\Gamma_1^{-1}+\Gamma_1
\biggl(\frac{B^*_{12}(t)}{(z-a_1(t))^2}+\frac{B^*_{22}(t)}
{(z-a_2(t))^2}+\sum_{i=1}^n\frac{B^*_{i1}(t)}{z-a_i(t)}\biggr)
\Gamma_1^{-1}.
$$

As the matrix $\Gamma_1(z,t)$ is holomorphically invertible (with
respect to $z$) in $\mathbb C$, one has
$$
B_{i1}(t)=\Gamma_1(a_i(t),t)\,B^*_{i1}(t)\Gamma^{-1}_1(a_i(t),t),
\qquad i=3,\ldots,n,
$$
and for $i=1,2$ one has
\begin{eqnarray*}
B_{i2}(t)&=&\Gamma_1(a_i(t),t)\,B^*_{i2}(t)\Gamma^{-1}_1(a_i(t),t),\\
B_{i1}(t)&=&\frac{\partial\Gamma_1}{\partial z}(a_i(t),t)\,
B^*_{i2}(t)\Gamma^{-1}_1(a_i(t),t)+\Gamma_1(a_i(t),t)\,B^*_{i1}(t)
\Gamma^{-1}_1(a_i(t),t)+ \\ & & +\Gamma_1(a_i(t),t)\,B^*_{i2}(t)
\frac{\partial\Gamma^{-1}_1}{\partial z}(a_i(t),t).
\end{eqnarray*}
Since
$$
\Gamma_1(z,t)=U'_0(t)^{-1}\Gamma'_1(z,t)=
\left(\begin{array}{c c}
\frac{f(t)}{u_1(t)} & -u_1(t) \\
\frac1{u_1(t)} & 0
\end{array}\right)\left(\begin{array}{c c}
1 & 0 \\
\frac{-z}{u_1(t)} & 1
\end{array}\right)=\left(\begin{array}{c c}
z+\frac{f(t)}{u_1(t)} & -u_1(t) \\
\frac1{u_1(t)} & 0 \end{array}\right)
$$
and the matrices $B^*_{ij}(t)$ are holomorphic near the point
$t=t^*$, one sees that the same holds for all the matrices
$(u_1(t))^2B_{ij}(t)$. {\hfill $\Box$}
\medskip

{\bf Remark.} Recall that the Painlev\'e III and V equations can
be described in terms of isomonodromic deformations satisfying
Theorem 3 (see details in \cite[Ch. 5, \S\S 4,5]{FIKN}): for $\rm
P_{III}$ one has $m=n=2$ and for $\rm P_{V}$ one has $m=1$, $n=3$.
If $t^*\in\Theta$ and $E|_{\overline{\mathbb C}\times \{t^*\}}
\cong\mathcal{O}(-k)\oplus\mathcal{O}(k)$, then the estimate
$2k\leqslant m+n-2$ holds \cite{Bol_Malek_Mit} when the monodromy
of a connection is irreducible. Thus $2k\leqslant2$ and hence
$k=1$ in the both cases.

\end{document}